\documentclass[12pt,a4paper]{article}
\usepackage{stmaryrd}
\usepackage{amsmath}
\usepackage{amsfonts}
\usepackage{amsmath,amsthm,dsfont,hyperref,latexsym}
\usepackage{graphicx}
\newcommand{\new}{\newcommand*}
\new{\rnew}{\renewcommand*}
\new{\newe}{\newenvironment*}
\new{\newt}{\newtheorem}
\new{\stl}{\setlength}
\new{\bea}{\begin{eqnarray}}
\new{\eea}{\end{eqnarray}}

\new{\be}{\begin{equation}}
\new{\ee}{\end{equation}}

\new{\bean}{\begin{eqnarray*}}
\new{\eean}{\end{eqnarray*}}
\new{\ba}{\begin{array}}
\new{\ea}{\end{array}}

\new{\no}{\nonumber}

\new{\bt}{\begin{theorem}}
\new{\et}{\end{theorem}}
\new{\bl}{\begin{lemma}}
\new{\el}{\end{lemma}}

\new{\bc}{\begin{corollary}}
\new{\ec}{\end{corollary}}
\new{\bp}{\begin{proof}\quad}
\new{\ep}{\end{proof}}

\rnew{\theequation}{\thesection.\arabic{equation}}
\rnew{\thefigure}{\thesection.\arabic{figure}}
\new{\sect}[1]{ \section{#1}
\setcounter{equation}{0} \setcounter{figure}{0} }
\newtheorem{lemma}{Lemma}[section]

\newtheorem{corollary}{Corollary}[section]
\newtheorem{theorem}{Theorem}[section]

\numberwithin{equation}{section} \numberwithin{table}{section}

\begin{document}

\centerline{\Large\bf Convex Bodies With Minimal Volume  } \vskip
10pt \centerline{\Large\bf Product in $\mathbb{R}^2$ --- A New Proof
} \vskip 25pt

\centerline{\large Lin Youjiang and Leng Gangsong} \vskip 15pt

\centerline{\footnotesize\sl$^{a}Department$ of Mathematics,
  Shanghai University, Shanghai 200444, P. R. China }
\centerline{\footnotesize linyoujiang@shu.edu.cn,
gleng@staff.shu.edu.cn}
$\;\;$\\

\footnote{

 \noindent{\footnotesize {2010 {\it Mathematics Subject Classification.}
Primary: 52A10, 52A40.}}

\noindent{\footnotesize\it Key words and phrases.} {\footnotesize
Convex body, Duality, Mahler Conjecture, Polytopes. }

\noindent{\footnotesize The authors would like to acknowledge the
support from the National Natural Science Foundation of China
(10971128), Shanghai Leading Academic Discipline Project (S30104).
}}



\noindent{\small {\bf Abstract.}  In this paper, a new proof of the
following result is given: The product of the volumes of an origin
symmetric convex bodies $K$ in $\mathbb{R}^2$ and of its polar body
is minimal if and only if $K$ is a parallelogram.
}\\


\section*{\normalsize\center 1. Introduction}
\setcounter{section}{1} \setcounter{equation}{0}
 $\;$

A well-known problem in the theory of convex sets is to find a lower
bound for the product of volumes $\mathcal {P}(K)=V(K)V(K^{\ast})$,
which is called {\it the volume-product of K }, where $K$ is an
$n$-dimensional origin symmetric convex body and $k^{\ast}$ is the
polar body of $K$ (see definition in Section 2).
Is it true that we always have
\begin{eqnarray}
\mathcal {P}(K)&\geq&\mathcal {P}(B_{\infty}^{n}),
\end{eqnarray}
where $B_{\infty}^{n}=\{x\in \mathbb{R}^n:~|x_i|\leq 1,~1\leq i\leq
n\}$?

For some particular classes of convex symmetric bodies in
$\mathbb{R}^n$, a sharper estimate for the lower bound of $\mathcal
{P}(K)$ has been obtained. If $K$ is the unit ball of a normed
$n$-dimensional space with a 1-unconditional basis, J. Saint-Raymond
\cite{SR81} proved that $\mathcal {P}(K)\geq 4^n/n!$; the equality
case, obtained for $1-\infty$ spaces, is discussed in \cite{Me86}
and \cite{Re87}. When $K$ is a zonoid it was proved in \cite{GMR88}
and \cite{Re86} that the same inequality holds, with equality if and
only if $K$ is an $n$-cube.

 In \cite{BM87}, J. Bourgain and V. D. Milman proved that there exist some $c>0$ such that for every $n$
and every convex body $K$ of $\mathbb{R}^n$,
$$\mathcal {P}(K)\geq c^n\mathcal {P}(B_2^n).$$ The best known
constant $c=\frac{\pi}{4}$ is due to Kuperberg \cite{Ku08}.

In \cite{Ma39}, K. Mahler proved (1.2) when $n=2$. There are several
other proofs of the two-dimensional result, see for example the
proof of M. Meyer, \cite{Me91}, but the question is still open even
in the three-dimensional case.

In this paper, we present a new proof about the problem when $n=2$,
which is different from the proof in \cite{Ma39} and \cite{Me91}.
Firstly, we prove that any origin symmetric polygon satisfies the
conjecture. Then, using the continuity of $\mathcal {P}(K)$ with
respect to the Hausdorff metric, we can easily prove that the
conjecture is also correct for any origin symmetric convex bodies in
$\mathbb{R}^2$. For the three-dimensional case, the conjecture maybe
can be solved by use of the same idea.

 Finally, let us mention the problem of giving an upper bound to $\mathcal {P}(K)$; it was proved by L. A. Santal\'{o} \cite{Sa49}: $P(K)\leq \mathcal {P}(B_2^n)$,
where $B_2^n$ is the $n$-dimensional Euclidean unit ball. In
\cite{LZ97}, \cite{Me89} and \cite{Pe85}, it was shown that the
equality holds only if $K$ is an ellipsoid.

\section*{\bf\normalsize\center 2. Notations and background materials}
\setcounter{section}{2}
 \setcounter{equation}{0}
$\;$

As usual, $S^{n-1}$ denotes the unit sphere, $B^n$ the unit ball
centered at the origin, $o$ the origin and $\|\cdot\|$ the norm in
Euclidean $n$-space $\mathbb{R}^n$. If $x$, $y\in \mathbb{R}^n$,
then $\langle x,y\rangle$ is the inner product of $x$ and $y$.

If $K$ is a set, $\partial K$ is its boundary, $int\;K$ is its
interior, and $conv~K$ denotes its convex hull. Let
$\mathbb{R}^n\backslash K$ denote the complement of $K$, i.e.,
$\mathbb{R}^n\backslash K=\{x\in \mathbb{R}^n: x\notin K\}.$ If $K$
is a $n$-dimensional convex subset of $\mathbb{R}^n$, then $V(k)$ is
its volume $V_n(K)$.

Let $\mathcal {K}^n$ denote the set of convex bodies (compact,
convex subsets with non-empty interiors) in $\mathbb{R}^n$. Let
$\mathcal {K}^n_o$ denote the subset of $\mathcal {K}^n$ that
contains the origin in its interior. Let $h(K,\cdot):
S^{n-1}\rightarrow \mathbb{R}$, denote the support function of $K\in
\mathcal {K}^n_o$; i.e.,

\begin{eqnarray}
h(K,u)= \max\{u\cdot x:~x\in K\}, u\in S^{n-1},
\end{eqnarray}
and let $\rho (K,\cdot): S^{n-1}\rightarrow \mathbb{R}$, denote the
radial function of $K\in \mathcal {K}^n_o$; i.e.,
\begin{eqnarray}
\rho(K,u)= \max\{\lambda\geq 0:~\lambda u\in K\}, u\in S^{n-1}.
\end{eqnarray}

A linear transformation (or affine transformation) of $\mathbb{R}^n$
is a map $\phi$ from $\mathbb{R}^n$ to itself such that $\phi x~=~A
x$ (or $\phi x~=~A x + t$, respectively), where $A$ is an $n \times
n$ matrix and $t\in \mathbb{R}^n$. By definition, for any
parallelograms centered at the origin $ABCD$ and
$A^{\prime}B^{\prime}C^{\prime}D^{\prime}$, there always is an
linear transformation $\mathcal {A}$ taking $ABCD$ to
$A^{\prime}B^{\prime}C^{\prime}D^{\prime}$.

Geometrically, an affine transformation in Euclidean space is one
that preserves:

(1). The collinearity relation between points; i.e., three points
which lie on a line continue to be collinear after the
transformation.

(2) Ratios of distances along a line; i.e., for distinct collinear
points $P_1$, $P_2$, $P_3$, the ratio $|P_2-P_1| / |P_3-P_2|$ is
preserved.

If $K\in {K}^n_o$, we define the polar body of $K$, $K^{\ast}$, by
$$K^{\ast}=\{ x\in \mathbb{R}^n:~ x\cdot y \leq 1~, \forall y\in K\}.$$

It is easy to verify that (see p.44 in \cite{Sc93})
\begin{eqnarray}
h(K^{\ast},u)=\frac{1}{\rho(K,u)}~~~~~~and~~~~~~~
\rho(K^{\ast},u)=\frac{1}{h(K,u)}
\end{eqnarray}

If $P$ is a polygon, i.e., $P= conv\{p_1,\cdots,p_m\}$, where $p_i$
$(i=1,\cdots,m)$ are vertices of polygon $P$. By the definition of
polar body, we have
\begin{eqnarray}
P^{\ast}&=&\{x\in \mathbb{R}^2: x\cdot p_1\leq 1,\cdots, x\cdot
p_m\leq 1\}
\nonumber\\
&=&\bigcap_{i=1}^{m}\{x\in \mathbb{R}^2: x\cdot p_i \leq 1\},
\end{eqnarray}
which implies that $P^{\ast}$ is the intersection of $m$ closed
half-planes with exterior normal vector $p_i$ and the distance of
straight line $\{x\in \mathbb{R}^2: x\cdot p_i= 1\}$ from
 the origin is $1/\|p_i\|$. Thus, if $P$ is an inscribed polygon in a unit circle,
then $P^{\ast}$ is polygon circumscribed around the unit circle. In
the proof of Lemma 3.3, we shall make use of these properties.

For $K$, $L\in \mathcal{K}^n$ the Hausdorff distance is defined by
\begin{eqnarray}
d(K,L)=\min\{\lambda\geq0:~K\subset L+\lambda B^n,~L\subset
K+\lambda B^n\},
\end{eqnarray}
which can be conveniently defined by (see p.53 in \cite{Sc93})
\begin{eqnarray}
d(K,L)=\max_{u\in S^{n-1}}|h(K,u)-L(K,u)|,
\end{eqnarray}
therefore, a sequence of convex bodies $K_i$ converges to $K$ if and
only if the sequence of support function $h(K_i, \cdot)$ converges
uniformly to $h(K,\cdot)$.

In $\mathcal {K}^n_o$, the convergence of convex bodies is
equivalent to the uniform convergence of their radial functions.
Because the conclusion will be used in the proof of Lemma 3.5, we
prove this conclusion (this proof is due to Professor Zhang Gaoyong
and we listened his lecture in Chongqing).

Let $K\in \mathcal {K}^n_o$. Define
\begin{eqnarray}
r_1=\max_{u\in S^{n-1}}\rho(K,u),
\end{eqnarray}
\begin{eqnarray}
r_0=\min_{u\in S^{n-1}}\rho(K,u).
\end{eqnarray}
It is easily seen that
\begin{eqnarray}
r_1=\max_{u\in S^{n-1}}h(K,u),
\end{eqnarray}
\begin{eqnarray}
r_0=\min_{u\in S^{n-1}}h(K,u).
\end{eqnarray}
\\
{\bf{Lemma 2.1.}} If $K\in \mathcal {K}^n_o$, then
\begin{eqnarray}
\rho(K+tB^n,u)\leq \rho(K,u)+\frac{r_1}{r_0}t,
\end{eqnarray}
\begin{eqnarray}
|u\cdot v(x)|\geq \frac{r_0}{r_1},
\end{eqnarray}
where $x=u\rho(K,u)\in \partial K.$\\
\\
{\bf{Proof.}} For $x\in \partial K$, let $x^{\prime}$ be the point
on $\partial (K+tB^n)$ and has the same direction as $x$. Let
$u=x/\|x\|=x^{\prime}/\|x^{\prime}\|$. Then
$$\rho(K+tB^n,u)-\rho(K,u)=\|x^{\prime}-x\|.$$
Since $K$ and $K+tB^n$ are parallel, the projection length of
$x^{\prime}-x$ onto the normal $v(x)$ is less than $t$,
$$\|x^{\prime}-x\|\leq \frac{t}{|u\cdot v(x)|}.$$
There is
\begin{eqnarray}
|u\cdot v(x)|&=&\frac{|x\cdot v(x)|}{\|x\|}
\nonumber\\
&=&\frac{h(K,v(x))}{\|x\|}
\nonumber\\
&\geq& \frac{r_0}{r_1}.
\end{eqnarray}
The desired inequalities follow. $\Box$
\\
\\
{\bf{Theorem 2.2.}} If a sequence of convex bodies $K_i\in \mathcal
{K}^n_0$ converges to $K\in \mathcal {K}^n_0$ in the Hausdorff
metric, then the sequence of radial functions $\rho(K_i,\cdot)$
converges to $\rho(K,\cdot)$ uniformly.
\\
{\bf{Proof.}} Assume that $d(K_i,K)<\varepsilon$. Then $K_i\subset
K+\varepsilon B^n$, and $K\subset K_i+\varepsilon B^n$. By Lemma
2.1,  (2.9) and (2.10)
$$\rho(K_i,\cdot)\leq \rho(K,\cdot)+\frac{r_1}{r_0}\varepsilon,$$
$$\rho(K,\cdot)\leq \rho(K_i,\cdot)+\frac{r_1+\varepsilon}{r_0-\varepsilon}\varepsilon.$$
When $\varepsilon<r_0/2$, we have
$$|\rho(K_i,\cdot)-\rho(K,\cdot)|\leq \frac{4r_1}{r_0}\varepsilon,$$
therefore the sequence of radial functions $\rho(K_i,\cdot)$
converges to $\rho(K,\cdot)$ uniformly. $\Box$

\section*{\bf\normalsize\center 3. Main result and its proof}
\setcounter{section}{3}
 \setcounter{equation}{0}
$\;$

 First, looking the following important theorem:\\

\noindent{\bf{Theorem 3.1.}} For any origin symmetric convex body
$K\subset \mathbb{R}^n$, $\mathcal {P}(K)$ is linear invariant, that
is, for every linear transformation
$A:~ \mathbb{R}^n\rightarrow \mathbb{R}^n$, we have $\mathcal {P}(AK)=\mathcal {P}(K)$. \\
 \\
{\bf{Proof.}} For any $u\in S^{n-1}$, we have
$$\rho((AK)^{\ast},u)=\frac{1}{h(AK,u)}=\frac{1}{h(K,A^{t}u)}=\rho(K^{\ast},A^tu)=\rho(A^{-t}K^{\ast},u).$$
Hence, $(AK)^{\ast}=A^{-t}K^{\ast}$,  therefore $$\mathcal
{P}(AK)=V(AK)V((AK)^{\ast})=V(AK)V(A^{-t}K^{\ast})$$$$=|A||A^{-t}|V(K)V(K^{\ast})=V(K)V(K^{\ast})=\mathcal
{P}(K).$$

$\Box$
\\

Because any parallelogram can been linear transformed into a unit
square, therefore their volume product is same (this value is equal
to 8). By the theorem above, we consider linear transformation of
origin symmetric polygon. We obtain the following theorem, which is
critical in our proof.
\\
\\
{\bf{Theorem 3.2.}} In $\mathbb{R}^2$, for any origin symmetric
polygon $P$, there exists a linear transformation $\mathcal
{A}:~P\rightarrow P^{\prime}$, where $P^{\prime}$ satisfies that
$P^{\prime}\subset B^2$ and there exist
three continuous vertices contained in $\partial B^2$.  \\
\\
{\bf{Proof}}. Since $P$ is origin symmetric polygon, its number of
sides is an even and corresponding two sides are parallel. Let
$A_1,\cdots,A_n, B_1,\cdots B_n$ denote all vertices of $P$. In
order to prove this theorem, we need three steps.

The first step,  transforming  parallelogram $A_1A_2B_1B_2 $ into
rectangular $A_1^{\prime}A_2^{\prime}B_1^{\prime}B_2^{\prime}$
inscribed in $B^2$. Now $P$ is transformed into $P_1$ (see (2) or
$(2)^{\prime}$ in Figure 3.1.1 and 3.1.2).
\begin{figure}[htb]
\centering
 \includegraphics[height=10cm]{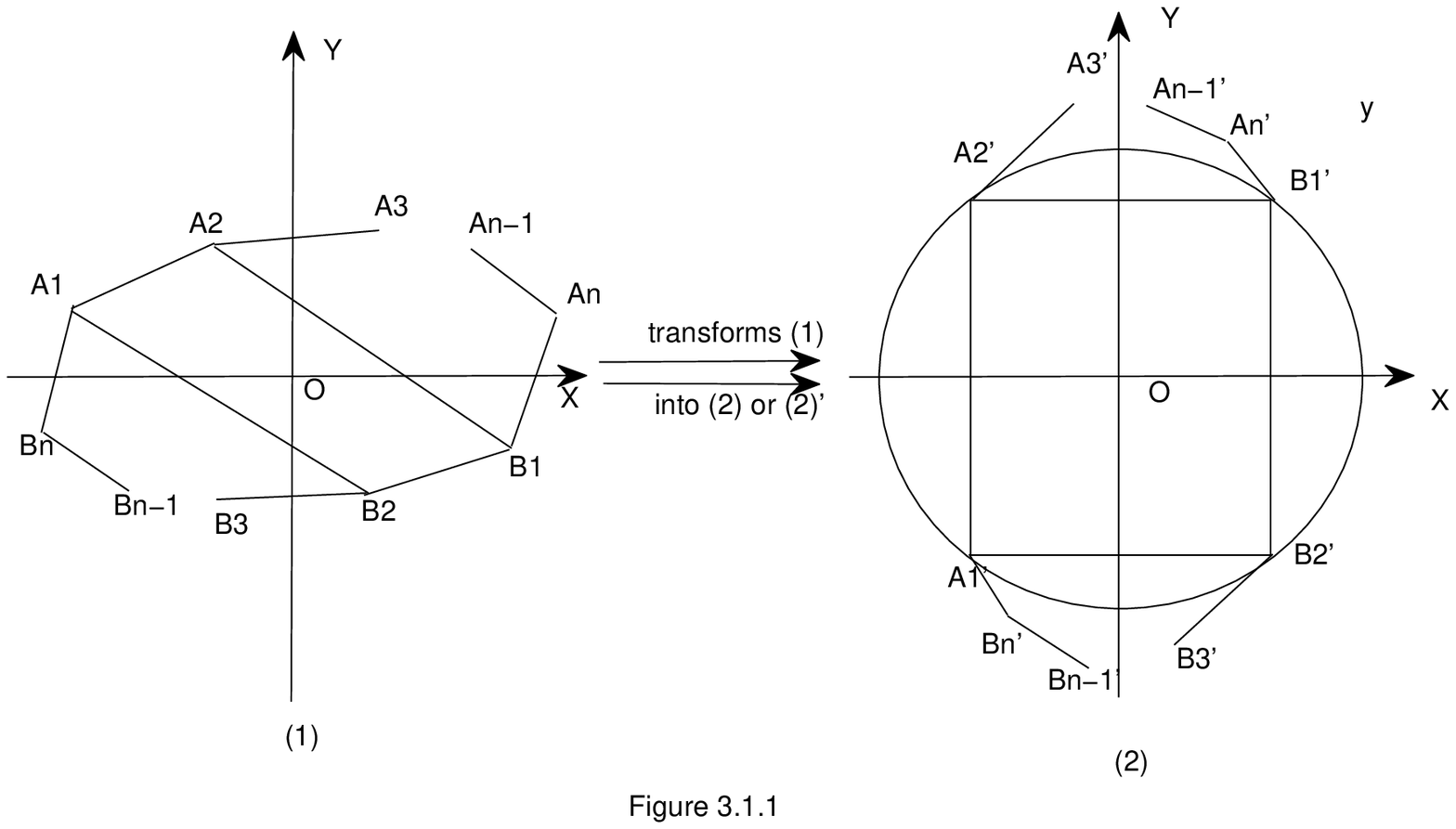}
\end{figure}
\begin{figure}[htb]
\centering
 \includegraphics[height=10cm]{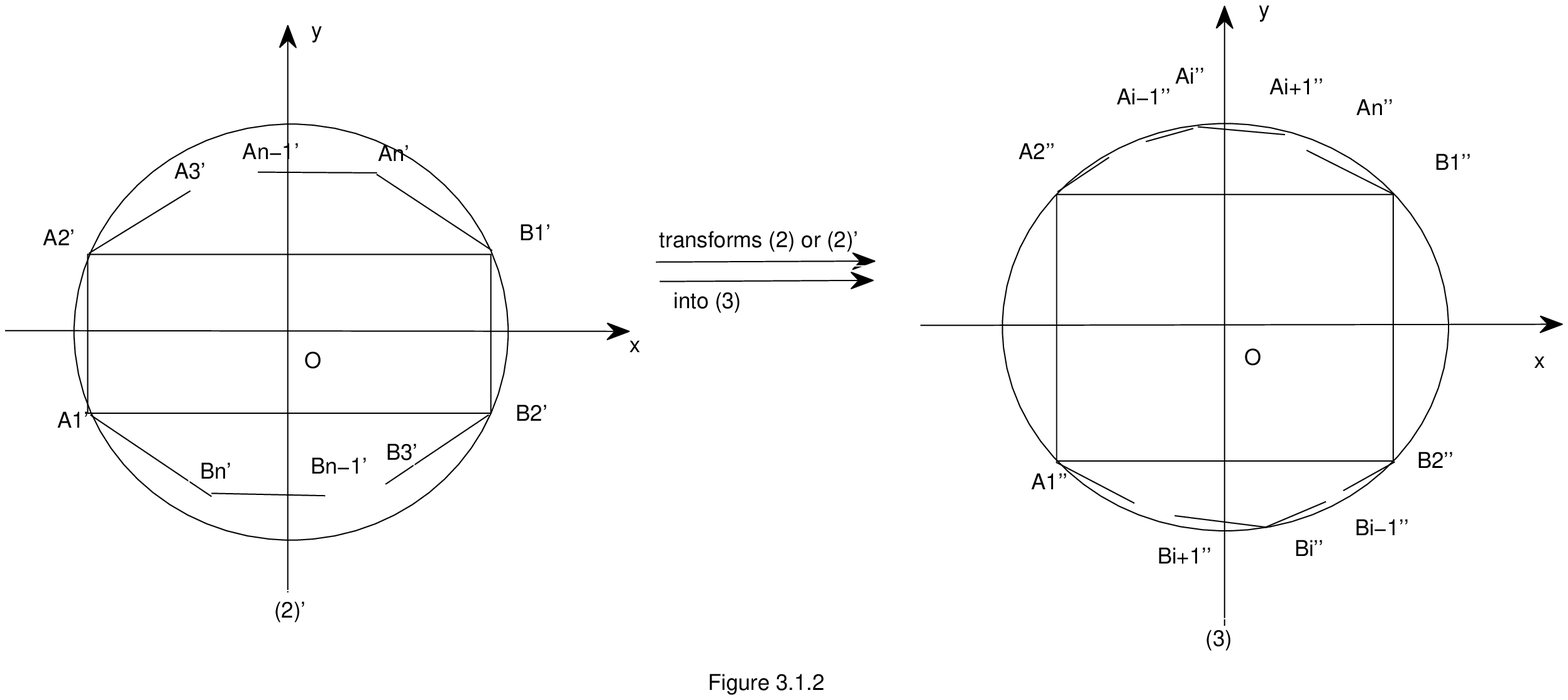}
\end{figure}

The second step, transforming $P_1$ into $P_2$ (see (3) in Figure
3.1.2). For polygon $P_1$, if there exist some vertices
 $$\{A_i^{\prime}:~i\in I\subset \{3,\cdots,n\}\}\subset \mathbb{R}^2\backslash B^2,$$
  there  exists a linear transformation $\mathcal {A}_1:P_1\rightarrow P_2$, which shortens segment
 $A_1^{\prime}A_2^{\prime}$ and $B_1^{\prime}B_2^{\prime}$ into $A_1^{\prime\prime}A_2^{\prime\prime}$ and
 $B_1^{\prime\prime}B_2^{\prime\prime}$, simultaneously makes some vertices
 $\{A_i^{\prime\prime}:i\in I_2\subset \{3,\cdots,n\}
 \}$ on boundary of $B^2$ and $P_2\subset B^2$.
If
 $$\{A_3^{\prime},\cdots,A_n^{\prime},B_3^{\prime},\cdots,B_n^{\prime}\}\subset int~B^2,$$
 then there exists a linear transformation $\mathcal {A}_1^{\prime}:P_1\rightarrow P_2$, which lengthens
 segments
 $A_1^{\prime}A_2^{\prime}$ and $B_1^{\prime}B_2^{\prime}$ into $A_1^{\prime\prime}A_2^{\prime\prime}$ and
 $B_1^{\prime\prime}B_2^{\prime\prime}$ respectively, simultaneously makes some vertices
 $\{A_i^{\prime\prime}:i\in I_1\subset \{3,\cdots,n\}
 \}$ on boundary of $B^2$ and $P_2\subset B^2$.  (see (3) in Figure 3.1.2).

The third step, transforming $P_2$ into $P_3$. If
$A_1^{\prime\prime},A_2^{\prime\prime},A_i^{\prime\prime}$ are three
continuous vertices contained in $\partial B^2$, then this theorem
has been proved; otherwise rotation transforming $P_2$ into
$P_3^{\prime}$, which satisfies that
$A_2^{\prime\prime}A_i^{\prime\prime}$ parallels x-axis (see (4) in
Figure 3.2). Then we transform $P_3^{\prime}$ into $P_3$,
lengthening segments $A_2^{\prime\prime}B_i^{\prime\prime}$ and
$A_i^{\prime\prime}B_2^{\prime\prime}$ into $A_2^{(3)}B_i^{(3)}$ and
$A_i^{(3)}B_2^{(3)}$ respectively, simultaneously making some
vertices
 $\{A_j^{(3)}:j\in I_3\subset \{3,\cdots,i-1\}
 \}$ on boundary of $B^2$ and $P_3\subset B^2$ (Since it is easy to prove that vertices
$\{A_{i+1}^{(3)},\cdots,A_n^{(3)},B_1^{(3)}\}$ are in the
 internal of $B^2$ ) (see (5) in Figure 3.2).

\begin{figure}[htb]
\centering
 \includegraphics[height=10cm]{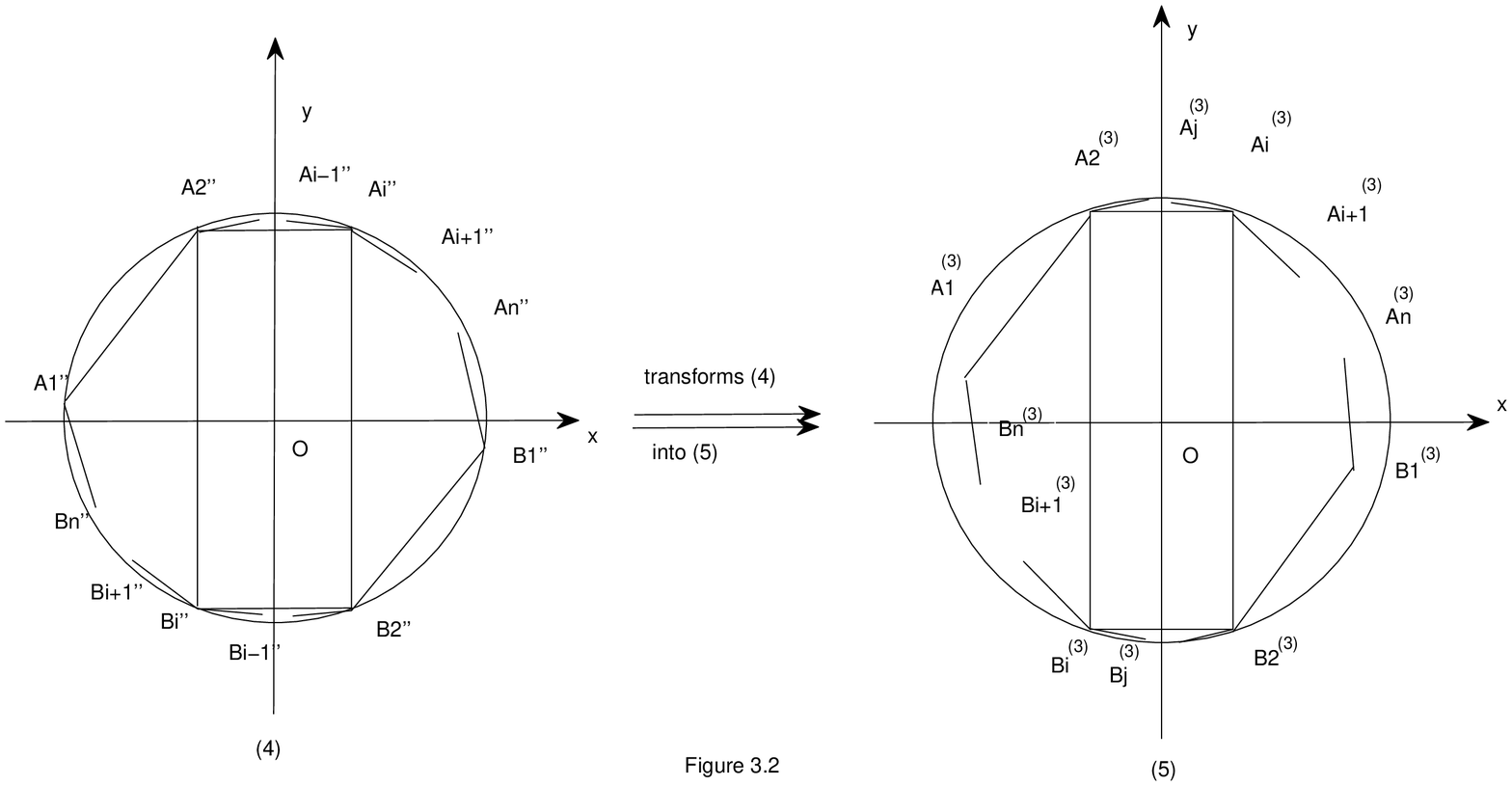}
\end{figure}

Repeating the third step finite times, we can get a polygon
$P^{\prime}$, in which there exist three continuous vertices
contained in $\partial B^2$, which completes the proof. $\Box$
\\

By above theorem, we consider the volume-product of polygon with
three continuous vertices in $\partial B^2$.
\\
{\bf{Lemma 3.3.}} Suppose that $P^{\prime}\subset B^2$ is an origin
symmetric polygon and  $A, C, B$ are three continuous vertices of
$P^{\prime}$ contained in $\partial  B^2$, then $\mathcal
{P}(P^{\prime\prime})\leq \mathcal {P}(P^{\prime})$, where
$P^{\prime\prime} $ is a new polygon from $P^{\prime}$ by deleting
vertices $C$ and $C^{\prime}$.\\
\\
{\bf{Proof.}} Suppose side $AB$ parallels X-axis (see Figure 3.3.),
straight lines $l$, $l_1$ and $l_2$ are  tangent lines to the unit
circle $B^2$ passing through points $C$, $A$ and $B$ respectively.
\begin{figure}[htb]
\centering
 \includegraphics[height=10cm]{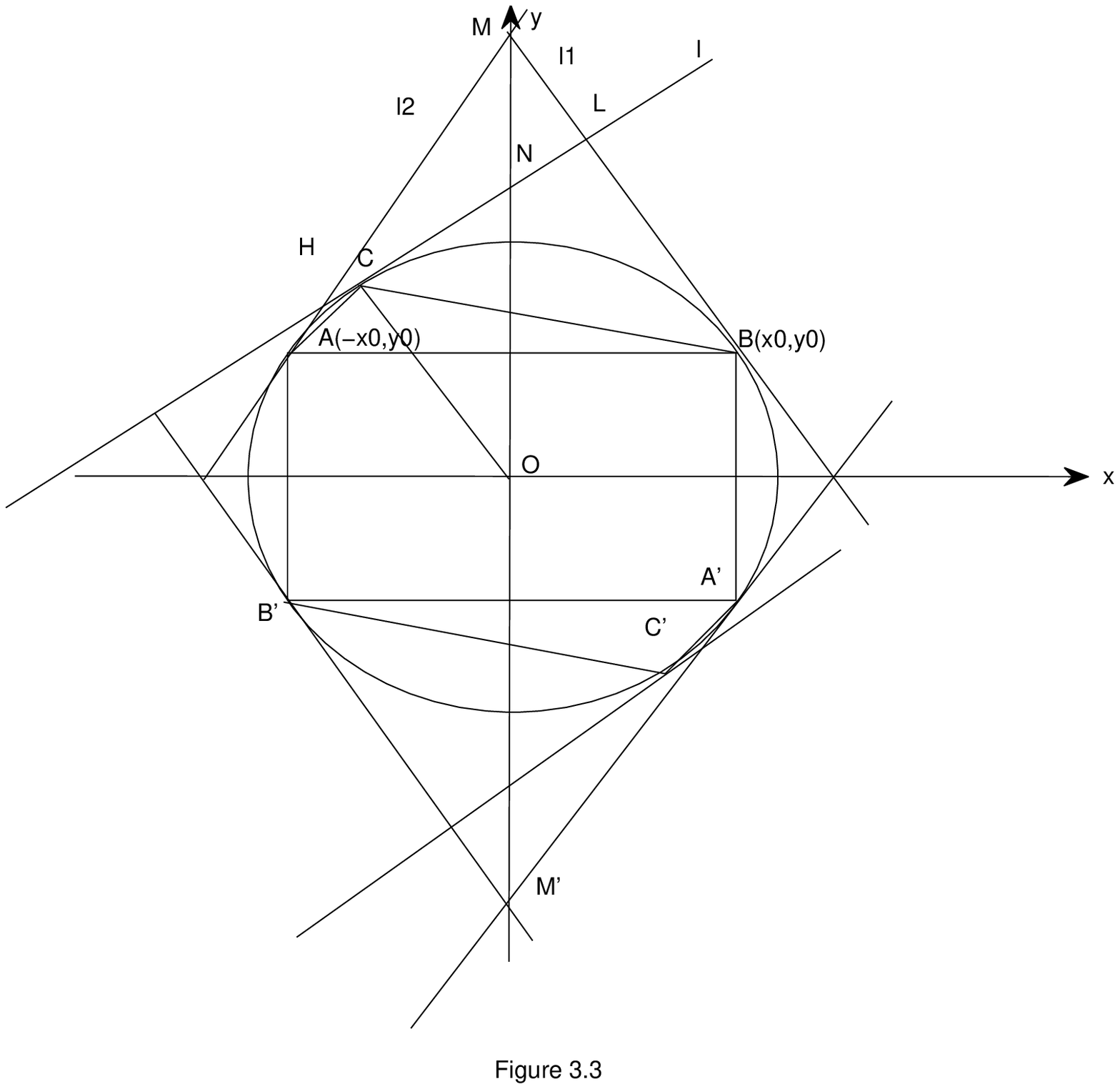}
\end{figure}
Let $A=(-x_0,y_0)$, then $B=(x_0,y_0)$. Let $\theta$ denote $\angle
xOC$. It is clear that $\pi/2\leq \theta\leq \pi-\arctan(y_0/x_0)$
when point $C$ is in third quadrant. We have the following equations
of straight lines:
$$l_1:~~y-y_0=\frac{x_0}{y_0}(x+x_0),$$
$$l_2:~~y-y_0=-\frac{x_0}{y_0}(x-x_0),$$
$$l:~~y-\sin \theta =-\frac{\cos \theta}{\sin \theta}(x-\cos \theta).$$
Let point $N$ denote the intersection of $l$ and Y-axis and point
$M$ denote the intersection of $l_1$ and Y-axis. We can easily get
$N(0,1/\sin \theta)$ and $M(0,1/y_0)$. In order to obtain the
abscissas of intersection of $l$ and $l_1$, $l_2$, we solve the
following equation systems:

\begin{equation} \label{eq:1}
\left\{ \begin{aligned}
         y-\sin\theta &= -\frac{\cos\theta}{\sin\theta}(x-\cos\theta) \\
                y-y_0  &=\frac{x_0}{y_0}(x+x_0)
                          \end{aligned} \right.
                          \end{equation}
and
\begin{equation} \label{eq:1}
\left\{ \begin{aligned}
         y-\sin\theta &= -\frac{\cos\theta}{\sin\theta}(x-\cos\theta) \\
                y-y_0  &=-\frac{x_0}{y_0}(x-x_0)
                          \end{aligned} \right.
                          \end{equation}

We can get abscissas of points $H$ and $L$:
$$x_1=\frac{y_0-\sin\theta}{y_0\cos \theta+x_0\sin \theta}$$
and $$x_2=\frac{y_0-\sin\theta}{y_0\cos\theta-x_0\sin\theta}.$$
Therefore we can obtain the area of $\triangle MHL$:
$$S_{\triangle MHL}=\frac{x_0}{y_0}\cdot\frac{\sin\theta-y_0}{\sin\theta+y_0}.$$

Let $V=V(P^{\prime\prime})$ and
$V^{0}=V({P^{\prime\prime}}^{\ast})$, where $P^{\prime\prime}$
denotes the new polygon from $P^{\prime}$ by deleting vertices $C$
and $C^{\prime}$, then $\mathcal {P}(P^{\prime})$ is a function
$f(\theta)$, where

\begin{eqnarray}
f(\theta)&=&\left(V+2x_0(\sin\theta-y_0)\right)\left(V^{0}-\frac{2x_0}{y_0}\cdot\frac{\sin\theta-y_0}{\sin\theta+y_0}\right)
\end{eqnarray}
and
$$\frac{\pi}{2}\leq \theta\leq
\pi-\arctan(\frac{y_0}{x_0}).$$

We have
\begin{eqnarray}
f^{\prime}(\theta)&=&2x_0\cos\theta\cdot\frac{(V^{0}y_0-2x_0)(\sin\theta+y_0)^2+2y_0(4x_0y_0-V)}{y_0(\sin\theta+y_0)^2}.
\end{eqnarray}
In (3.6), since $\cos\theta\leq 0$ and $y_0(\sin\theta+y_0)^2\geq
0$, in order to prove $f^{\prime}(\theta)\leq 0$, let
$t=\sin\theta$, we just need to prove $g(t)\geq 0$, where

\begin{eqnarray}
g(t)&=&(V^{0}y_0-2x_0)(t+y_0)^2+2y_0(4x_0y_0-V),~~~~t\in [y_0,1].
\end{eqnarray}

In order to prove $g(t)\geq 0$, we just need to prove that
$V^{0}y_0-2x_0>0$ and $g(y_0)\geq 0$.

Because of $V({P^{\prime\prime}}^{\ast})\geq V(conv\{
A,M,B,A^{\prime},M^{\prime},B^{\prime}\})$,
$$V^{0}\geq 4x_0y_0+2x_0(\frac{1}{y_0}-y_0),$$ therefore,
\begin{eqnarray}
V^{0}y_0-2x_0&\geq&\left(4x_0y_0+2x_0(\frac{1}{y_0}-y_0)\right)y_0-2x_0
\nonumber\\
&=&2x_0y_0^2
\nonumber\\
&>&0,
\end{eqnarray}
and therefore function $g(t)$ is a parabola opening upward. Thence,
when $t\in [y_0,1]$, quadratic function $g(t)$ is increasing, thus
we just need to proof
\begin{eqnarray}
g(y_0)&=&2y_0(2V^0y_0^2-V)\geq 0.
\end{eqnarray}

 Let $\mathcal {D}$ denote the area of circular segment enclosed by
arc $\widehat{BA^{\prime}}$ and chord $\overline{BA^{\prime}}$, then
\begin{eqnarray}
V^{0}&\geq&4x_0y_0+2x_0(\frac{1}{y_0}-y_0)+2\mathcal {D}
\end{eqnarray}
and
\begin{eqnarray}
V&\leq&4x_0y_0+2\mathcal {D}.
\end{eqnarray}

In order to prove (3.9), we just need to prove
\begin{eqnarray}
2\left(4x_0y_0+2x_0(\frac{1}{y_0}-y_0)+2\mathcal
{D}\right)y_0^2&\geq&4x_0y_0+2\mathcal {D},
\end{eqnarray}
which equivalent to
\begin{eqnarray}
2x_0y_0^3&\geq&\mathcal {D}(1-2y_0^2).
\end{eqnarray}
And because
\begin{eqnarray}
\mathcal {D}&\leq&(1-x_0)\cdot 2y_0,
\end{eqnarray}
hence, we just need to prove
\begin{eqnarray}
x_0y_0^3&\geq&y_0(1-x_0)(1-2y_0^2),
\end{eqnarray}
which equivalent to
\begin{eqnarray}
x_0^3-2x_0^2+1&\geq&0,
\end{eqnarray}
which is clearly correct.

Summary, we get $f^{\prime}(\theta)\leq 0$ when $\theta\in
[\pi/2,\pi-\arctan(y_0/x_0)]$, hence when
$\theta=\pi-\arctan(y_0/x_0)$, which implies that point $C$
coincides with point $A$, function $f(\theta)$ obtain minimal
function value, therefore $\mathcal
{P}(P^{\prime\prime})\leq\mathcal {P}(P^{\prime})$.
 $\Box$\\

Making use of Lemma 3.3, we can obtain the following conclusion.
\\
\\
{\bf{Theorem 3.4.}}  If $P\subset \mathbb{R}^2$ is an origin
symmetric polygon, then $\mathcal {P}(P)\geq
\mathcal {P}(S)$, where $S$ is square.\\
\\
{\bf{Proof.}} By Theorem 3.2, Lemma 3.3 and linear invariance of
$\mathcal {P}(P)$, if the number of sides of polygon $P$ is $2n$,
there exists a polygon $P_1$ with $2(n-1)$ sides satisfying
$\mathcal {P}(P_1)\leq \mathcal {P}(P)$. Repeating this process
$n-2$ times, we can obtain a square $S$ satisfying
$\mathcal{P}(P)\geq \mathcal{P}(S)$. $\Box$
\\

\indent In order to obtain the main result in the paper, we first
prove the following lemma.
\\
\\
{\bf{Lemma 3.5.}} The volume product $\mathcal {P}(K)$ is continuous
under the Hausdorff metric.\\
\\
{\bf{Proof.}} Let
$$\lim_{i\rightarrow\infty}K_i=K.$$
By Theorem 2.2, the sequence of radial function $\rho(K_i,\cdot)$
converges to $\rho(K,\cdot)$ uniformly, therefore the reciprocal of
radial function $1/\rho(K_i,\cdot)$ converges to $1/\rho(K,\cdot)$
uniformly. Since
\begin{eqnarray}
d(K_i^{\ast},K^{\ast})&=& \max_{u\in
S^{n-1}}|h(K^{\ast}_i,u)-h(K^{\ast},u)|
\nonumber\\
&=&\max _{u\in
S^{n-1}}\left|\frac{1}{\rho(K_i,u)}-\frac{1}{\rho(K,u)}\right|,
\end{eqnarray}

we have
\begin{eqnarray}
\lim_{i\rightarrow \infty}K^{\ast}_i&=&K^{\ast}.
\end{eqnarray}
By continuity of the volume function $V(\cdot)$ under the Hausdorff
metric, we have
\begin{eqnarray}
\mathcal {P}(K)&=&V(K)V(K^{\ast})
\nonumber\\
&=&\lim_{i\rightarrow\infty}V(K_i)\lim_{i\rightarrow\infty}V(K_i^{\ast})
\nonumber\\
&=&\lim_{i\rightarrow\infty}V(K_i)V(K_i^{\ast})
\nonumber\\
&=&\lim_{i\rightarrow\infty}\mathcal {P}(K_i).
\end{eqnarray}
$\Box$\\
{\bf{Theorem 3.6.}} If $K\subset \mathbb{R}^2$ is an origin symmetric convex body and $S\subset \mathbb{R}^2$ is a square, then $\mathcal {P}(K)\geq \mathcal {P}(S)$.\\
\\
{\bf{Proof.}} For any origin symmetric convex body $K\subset
\mathbb{R}^2$, there exists a sequence of origin symmetric polytopes
$\{P_i\}$ converging to $K$ under the Hausdorff metric. By Theorem
3.4 and Lemma 3.5, we have
\begin{eqnarray}
\mathcal {P}(K)=\lim_{n\rightarrow \infty}\mathcal {P}(P_i)\geq
\mathcal {P}(S).
\end{eqnarray}
$\Box$\\


\end{document}